\documentclass{article}
\usepackage[T2A]{fontenc}
\usepackage[english]{babel}
\usepackage[tbtags]{amsmath}
\usepackage{amsfonts,amssymb,mathrsfs,amscd,amsmath,comment}

\begin{document}

\author{B.\,S.~Safin\footnote{Nizhnij Novgorod, Russia; E-mail: {\tt boris-nn@yandex.ru}}}
\title {Description of Euler bricks using~Fibonacci's~identity}

\maketitle
\begin{abstract}
We show how the Fibonacci's identity is used to obtain Euler bricks.
Also, we put forward the relation between Fibonacci's identity and Euler's formula,
 $(X=n^6-15n^4+15n^2-1$, $Y=6n^5-20n^3+6n$, $Z=8n^5-8n)$,
which provides the description of Euler's bricks with noninteger spatial diagonal.
Finally, we establish a~relation between the Euler bricks with integer and noninteger spatial diagonals.
\smallskip

\end{abstract}

A rectangular parallelepiped for which all elements (sides and diagonals) are integers, except for one element, is called
an Euler brick. A~perfect cuboid is, by definition, a~rectangular parallelepiped whose elements are all integers.
A~brick is a~rectangular parallelepiped with at least two noninteger elements.
A~Pythagorean triples $(a,b,c)$ is called simplest if $\operatorname{g.c.d.}\, (a,b,c) =1$;
the corresponding triangle is also called simplest.444
The rule to obtain Euler bricks will be conventionally called the multiplication of Pythagorean triangles.

The book by Davenport \cite {Dav} (page 116 of the Russian translation) contains the following form of Fibonacci's identity
\begin{gather}
  \begin{matrix}
  a^2 +b^2=c^2,\\
  d^2+e^2=f^2
  \end{matrix}\qquad \qquad
 \begin{matrix}
  ae +db ,\\
  ae -db,\\
  ad+be,\\
 ad-be,
  \end{matrix} \nonumber
  \\
\begin{aligned}
  (a^2+b^2)(d^2+e^2) &= (ae+db)^2+(ad-be)^2 =\\
                     & = (ae-db)^2+(ad+be)^2=
         (cf)^2 ; 
\end{aligned}
\end{gather}
this identity is valid, in particular, for legs $a,b,d,e$ of right Pythagorean triangles,
$a^2 +b^2=c^2$, $d^2+e^2=f^2$.

Assume that the greatest common divisor of Pythagorean triples  $(a,b,c)$ and $(e,d,f)$ is~1.

Let $a, d$ be even numbers. Assume that $c, f$ can be written as a~single sum of two.
These assumptions are not essentials and are used to keep the resulting makes the analysis significantly simpler.

In principle, this is the modification of the Fibonacci identity in geometric interpretation (note that it can also
be extended to a~larger number of terms; i.e., to multiply an arbitrary number of Pythagorean triangles). Let us multiply
the second and third triangles.

What can we derive from this construction (here we multiply two triangles):

1. Firstly, we obtain two pairs of new Pythagorean triangles:
$$
\begin{matrix}
 (ae+db)^2+(ad-be)^2   = (cf)^2, \\
  (fa)^2 +(fb)^2=(fc)^2
\end{matrix}\qquad \qquad
\begin{matrix}
(ae-db)^2+(ad+be)^2= (cf)^2, \\
 (cd)^2+(ce)^2=(cf)^2
\end{matrix}
$$
In these pairs, the triangles $(fa)^2 +(fb)^2=(fc)^2$ and $(cd)^2+(ce)^2=(cf)^2$ are interchangeable.
The principle upon which they are formed is unknown, and so one has invoke simple enumeration.

Note that the triangles
$$
 (ae+db)^2+(ad-be)^2 = (cf)^2\qquad \text{and}\qquad (ae-db)^2+(ad+be)^2= (cf)^2
$$
are known to be simplest (see~\cite{Se}).

If we assume that $cf$ is a~spatial diagonal of a Euler brick, then combining these triangles in pairs or triples
we may obtain the set of Euler bricks, most of which would contain several noninteger elements.
For a~while, we shall be concerned with bricks having one or two nonzero elements. Also, triangles will be taken in pairs,
because using a~pair of triangle with the same diagonal it already possible to create a~Euler brick, in which 5~elements are known,
 and 2~remaining ones can be determined.

2. Secondly, this construction yields four more Euler bricks with noninteger spatial diagonal.
In other words, we find four triples of numbers that determine four Euler bricks with noninteger spatial diagonal.
The general rule is as follows: we need to examine 4~triples $ (ae)$, $(db)$, $(ad)$, $(ae)$, $(db)$, $(be)$, $(ad)$, $(be)$, $(ae)$, $(ad)$, $(be)$, $(db)$. These triples will have one noninteger element (the spatial diagonal), provided that the square root of the numbers  $( ae) ^2+(db)^2$ or ($ad)^2+(be)^2$
is an integer number.

If the square root of the number $( ae) ^2+(db)^2 $ is integer, then we obtain the following triples:
$(ae)$, $(db)$, $(ad)$ and $(ae)$, $(db)$,~$(be)$. Otherwise, we need to examine the numbers $(ad)$, $(be)$, $(ae)$ and $(ad)$, $(be)$,~$(db)$.
Since the Fibonacci identity follows from the law of multiplication of conjugate complex numbers, the above rule can be derived
from the representation of a~complex number as a~real matrix of the form $\begin{pmatrix} x & y \\ -y & x\end{pmatrix}$
with usual matrix addition and multiplication.

Thus, by multiplying two Pythagorean triangles (whose hypotenuses can be uniquely expressed as the sum of
two squares), we obtain a~set of Euler bricks with integer spatial diagonal and 4~bricks with noninteger spatial diagonal.
In the present note, we shall ne concerned not with all resulting bricks, but rather with those that
have one noninteger element.

We note at once that multiplication of triangles not always produces all integers
(of course, with the exception of the spatial (lateral) diagonal or side), but it always results in Euler
bricks (however, no proof of this fact is available). The search of triangles that produce integer solutions (of course, except for one element) is a separate problem, which was solved in partial already by Euler. This matter will be address below, but now let us consider a~classical example of a brick with integer spatial diagonal to explain the crux of the matter.

Consider a triple of numbers that generates a~Euler brick with integer spatial diagonal: 104, 153, 672,
the spatial diagonal is as follows: $697=17\times 41$. Multiplying the triangles
$$
  \begin{matrix}
   8^2 + 15^2 =17^2 ,\hfill \\
 40 ^2 + 9 ^2  =41^2\hfill {}
 \end{matrix}\qquad
  \begin{matrix}
 8\times 9 +40\times 15 =72 +600=672, \hfill \\
 8\times 9 - 40\times 15 =72- 600=-528,\hfill  \\
 8\times 40+ 15\times 9  =320+135=455, \hfill \\
 8\times 40-15\times 9    =320-135=185,\hfill {}
\end{matrix}
$$
we obtain two triangles $672^2+185^2=528^2+455^2=697^2$, the second pair is obtained by multiplying the triangle $8^2 + 15^2 =17^2$ by 41, and the triangle $40 ^2 + 9 ^2  =41^2$, by 17. That is,
\begin{gather*}
 (41\times 8)^2+(41\times 15)^2=(41\times 17)^2 =328^2+615^2=697^2, \\
(17\times 40)^2+(17\times 9)^2=(17\times 41)^2=680^2+153^2=697^2.
\end{gather*}
Thus, combining in pairs, we obtain 4~Euler bricks with integer spatial diagonal:
\begin{align*}
&672^2+185^2=697^2,& \quad            &528^2+455^2=697^2, \\
&680^2+153^2=697^2,& \quad            & 328^2+615^2=697^2, \\
\noalign{\medskip}
&672^2+185^2=697^2,& \quad            &528^2+455^2=697^2, \\
&328^2+615^2=697^2,& \quad            & 680^2+153^2=697^2.
\end{align*}
Consider the first two pairs (the remaining two pairs are devoid of interest, because they have more
than two noninteger elements). For these pairs, we know five elements of Euler bricks; it remains to find the two missing ones:
\begin{align*}
&680^2-672^2=104^2,&\qquad                 &528^2-328^2~413{,}8^2, \\
&697^2-104^2=(689{,} 1)^2 ,& \qquad        &697^2 -413{,}8^2=(561{,} 1)^2
\end{align*}

Thus, we have obtained two Euler bricks, of which one has noninteger lateral diagonal,
and the other, a~noninteger side and noninteger lateral diagonal. For the first brick we have:
104, 672, 153 (the sides), 680, 185, (689, 1) (the lateral diagonals), and 697 (the spatial diagonal).

The second brick (2 noninteger elements) has two sides (413, 1), 328, 455, the lateral diagonals
528, 615, (561, 1), and the spatial diagonal 697.

The same construction is also capable of producing 4 more bricks with noninteger spatial diagonals.
Their sides are (72, 600, 320), (72, 600,~135), (320, 135,~72), (320,
135, 600).
In this particular example these bricks also have two
or more noninteger elements (the spatial diagonal and the lateral diagonal), but this is not relevant, because the
following example may be used to illustrate how all the bricks with noninteger spatial diagonal may be obtained.

Consider a~classical example of a brick with noninteger spatial diagonal (here we shall be concerned only with
bricks with one noninteger element). The sides are: 44, 117, 240. We obtain three
Pythagorean triangles:
\begin{gather*}
 44^2+117^2=125^2, \\
 44^2+240^2=244^2  =4^2(11^2+60^2=61^2), \\
 117^2+240^2=267^2    =3^2(39^2+80^2=89^2)
\end{gather*}
Multiplying the two last triangles, $11^2+60^2=61^2$, $39^2+80^2=89^2    $,
\begin{align*}
 &11\times 80+39\times 60=880+2340=3220, \\
  &            11\times 80-39\times 60=880-2340=-1460, \\
  &     11\times 39+60\times 80=429+4800=5229, \\
  &11\times 39-60\times 80=429-4800=-4371,
\end{align*}
we obtain two triples of numbers, which describe two Euler bricks with noninteger spatial diagonals: $880$, $2340$, $429$, and 880, 2340, 4800.
Cancelling the second triple by 20, we arrive at two complete triples: 44, 117,~240 and 429, 880,~2340. Figuratively speaking, each triple consists of three numbers, of which two are multiplied legs of different parity, plus two multiplied even legs for one triple and two multiplied odd legs for the other triple. In this example, we also obtain two bricks with integer diagonal: $3220^2+4371^2=1460^2+5229^2=5429^2$; however, these will not be entered into here, because, in my opinion, the number of noninteger elements is greater than two.

Consider another example: the triple: 140, 480,~693.

We have three equations:
\begin{align*}
&480^2+693^2=843^2 =3^2(160^2+231^2=281^2), \\
&480^2+140^2=500^2 =20^2(24^2+7^2=25^2), \\
&140^2+693^2=707^2 =7^2(20^2+99^2=101^2)
\end{align*}

Multiplying pairwisely the three last equations, we obtain the following triples:
5544, 1120, 3840 (cancelling by 8), 693, 140, 480, and 5544, 1120,~1617, of which the first triple is obtained from the first two equations.

Multiplying the second and third equations, we have the triples: 693, 480, 140 and 480,
693, 2376, which after cancelling by~3, gives 792, 231,~160.

Multiplying the firth and third equations we have: $4620$ $15840$, $3200$ or
(after cancelling by 20), 231, 792, 160, and 4620, 15840, 22869.

Having the notable Euler's formula at our disposal (see, e.g., \cite{Ostrik}), it is an easy matter to find the triples with
noninteger spatial diagonal (for example, $X=n^6-15n^4+15n^2-1$, $Y=6n^5-20n^3+6n$, $Z=8n^5-8n$). For $n=2$, the problem was already solved above.

For $n=3$ this gives even numbers; after cancelling, we arrive at $n=2$, while for $n=4$ we obtain
495, 4888, 8160. Using these numbers, we find the triangles:
\begin{align*}
&495^2 +4888^2 =4913^2 , \\
&495^2 +8160^2 =8175^2 ,&\quad &15^2(33^2 +544^2 =545^2), \hfill \\
&4888^2 +8160^2 =9512^2,&\quad  &8^2(611^2 +1020^2 =1189^2)\hfill{}
\end{align*}
Multiplying the triangles
\begin{align*}
& 33^2+544^2=545^2, \\
& 611^2+1020^2=1189^2 ,
\end{align*}
we arrive at the following two triples: 332384, 33660, 554880 and 332384, 33660, 20163. Cancelling the first triple by~68, we find two complete triples: 495, 8160, 4888 and    332384, 33660, 20163.

However, the problem of finding bricks with integer spatial diagonal involves much more great difficulties; but this is a~separate problem.

\smallskip

Consider one more example with integer spatial diagonal. Let us multiply pairwisely three Pythagorean triangles:
\begin{align*}
&3^2+4^2=5^2, \\
& 5^2+12^2=13^2, \\
&8^2+15^2=17^2
\end{align*}
Omitting the calculations, we shall state the final result: there are four new triangles:
\begin{equation}
\begin{aligned}
&1073^2+264^2=1105^2, \\
& 943^2+576^2=1105^2,\\
&817^2+744^2=1105^2,\\
&47^2+1104^2=1105^2
\end{aligned} 
\end{equation}
Since the legs and hypotenuses of these triangles have no common divisors, we obtain 10 more triangles, which have common divisors. For brevity, we shall not write them down, but point out that they are obtained according to case~1 by multiplication of the corresponding hypotenuses by the corresponding legs and hypotenuses.

Let us write down only one equation, which is essential for further explanation:
$$
520^2+975^2=1105^2   ;
$$
this equation is obtained from the equation $40^2+75^2=85^2 $ by multiplication by~13; the latter equation follows, in turn, from the equation $ 8^2+15^2=17^2 $ by multiplication by~5. Thus, consider two pairs of equations:
\begin{align*}
&1073^2+264^2=1105^2,&\qquad &943^2+576^2=1105^2, \\
&  520^2+975^2=1105^2,& \qquad & 520^2+975^2=1105^2    .
\end{align*}

From the first pair of equations we obtain: the sides (448, 264, 975), the lateral diagonals  1073, 520, (1010, 1), and, of course, the spatial diagonal 1105.

From the second pair we obtain: the sides 520, 576, (786, 7), the lateral diagonals 943, 776, 975, and the spatial diagonal 1105.

Consequently, it may be conjectured that the cases when bricks have a~noninteger lateral diagonal or one noninteger side may be combined into one case. Then there are only two different cases in the brick problem:

\smallskip

1. Noninteger spatial diagonal.

 2. Integer spatial diagonal.

\smallskip

Keeping all the assumptions for the moment, the logic of the proof of existence or nonexistence of a~brick amounts
to showing that it is impossible to obtain an integer spatial diagonal in the first case or and that a~side or
a~lateral diagonal will prove to be noninteger.

\smallskip

This immediately suggests the following questions:

1. Is the present construction capable of describing all the Euler bricks?

2. How to prove or disprove this.

3. We multiply arbitrary Pythagorean triangles---may be this is the desired proof that all the triangles are described?

4. Multiplication of three and more triangles produces four or more simplest Pythagorean triangles---how
to prove that they are not elements of one brick?

\medskip

Assume we are given the system of equations describing the perfect cuboid:
\begin{equation}\label{eq3} 
\begin{aligned}
&m^2+n^2=p^2, \\
&m^2+k^2=q^2  , \\
&n^2+k^2=z^2, \\
&n^2+q^2=g^2,
\end{aligned}
\end{equation}
where $(m,n,k)$ are legs of the cuboid, ($p,q,z$) are lateral diagonals, (g) is the spatial diagonal,
$m,n$ are even and $k$~is odd (see~\cite{Lee}).
\smallskip

For any right parallelepiped the following condition is satisfied:
$$
  m^2+n^2+k^2=g^2.
$$

Let $x_1y_1$ be two numbers generating the equation $m^2+n^2=p^2$, $(m=2x_1 y_1$, $n=x_1^2-y_1^2$, $p=x_1^2+y_1^2)$,
and let $x_2 y_2$ be two numbers generating the equation $ m^2+k^2=q^2$, $(m=2x_2 y_2$, $k=x_2^2-y_2^2$, $q=x^2_2+y_2^2)$.

Pythagorean triples can be found from the formulas: $(2xy,x^2-y^2,x^2+y^2)$, where $x$ and $y$ are integers.
Now the equation $m^2+n^2+k^2=g^2 $ can be rewritten as follows. We have $2x_1y_1=2x_2y_2$ and hence
$$
  4x^2_1y^2_1=4x^2_2y^2_2+(x^2_1-y^2_1)^2+(x^2_2-y^2_2)^2=g^2.
$$
Expanding, we obtain the following equation \cite{Lee}:
$$
   4x^2_1y^2_1=4x^2_2y^2_2+x_1^4-2x^2_1y^2_1+y_1^4+x_2^4-2x^2_2y^2_2+y_2^4=g^2=x_1^4+y_1^4+x_2^4+y_2^4.
$$
We have  $4x^2_1 y^2_1=4x^2_2y^2_2$, and to the expression $x_1^4+y_1^4+x_2^4+y_2^4$ we may add and subtract the expression $4x^2_1y^2_1=4x^2_2y^2_2$,
thereby obtaining the equation:
\begin{equation}\label{eq4} 
(x_1^2+y_1^2)^2+(x_2^2-y_2^2)^2=(x_1^2-y_1^2)^2+(x_2^2+y_2^2)^2=g^2            . 
\end{equation}
Recalling the Fibonacci identity
\begin{equation}\label{eq1n} 
(a^2+b^2)(d^2+e^2)=(ae+db)^2+(ad-be)^2=(ae-db)^2+(ad+be)^2=(cf)^2 
\end{equation}
and assuming that
\begin{equation}\label{eq5} 
ae=x_1^2,\quad db=y_1^2, \quad ad=x_2^2, \quad be=y_2^2,
\end{equation}
then the expression
\begin{equation}\label{eq4n} 
 (x_1^2+y_1^2)^2+(x_2^2-y_2^2)^2=(x_1^2-y_1^2)^2+(x_2^2+y_2^2)^2=g^2                      
\end{equation}
is also the Fibonacci identity. Hence, we may put forward the above scheme of obtaining Euler bricks; that~is:
$$
  \begin{aligned}
&ae +db  =  x_1^2+y_1^2, \\
&ae -db  =  x_1^2-y_1^2, \\
&ad +be  =  x_2^2+y_2^2, \\
&ad- be  =  x_2^2-y_2^2.
  \end{aligned}
$$
Since we assume that $ae $, $db$, $ad$, $be$ are legs of Pythagorean triangles (which form two Euler bricks,
though with noninteger spatial diagonal), hence the expressions
$(ae)^2+(db)^2$, $(ad)^2+(be)^2$, $(ae)^2+ (ad)^2$,  $(ae)^2 +(be)^2$ should be squares.

Consequently, it follows that the equations $ x_1^4+y_1^4=\text{square}$ or $x_2^4+y_2^4=\text{square}$
and other obeys ($x_1^4+ x_2^4$, $x_1^4+ y_2^4)$ may have a~solution in integers
(since $ae=x_1^2$, $db=y_1^2$, $ad=x_2^2$, $be=y_2^2$), which is impossible.

\smallskip

{\bf Corollary.} \textit{The equation $x^4+y^4+z^4=w^2$ is not solvable, provided that any of the integers $x, y, z$ is a~product
of legs of {\rm 2}~Pythagorean triangles.}

\smallskip

For the case involving an integer spatial diagonal, the following conditions must be satisfied:
$$
(ae+db)^2+(ad-be)^2=(ae-db)^2+(ad+be)^2= (cf)^2,
$$
because $ae=x_1^2$, $db=y_1^2$, $ad=x_2^2$, $be=y_2^2$, and so we arrive at equation~\eqref{eq4}.
It was \textit{ab initio} assumed that the triangles $  a^2 +b^2=c^2$,  $d^2+e^2=f^2$ are primitive,
hence the legs $(a,b)$, $(d,e)$ are mutually prime. According to the theorem on the product of two mutually
primes, if this prime is a~square, then the mutually prime numbers are also squares. Hence, $ a=r^2$, $b=s^2$, $d=u^2$, $e=w^2$.
Consequently, two Pythagorean triangles that satisfy the Fibonacci identity can be rewritten as follows:
$$
  (a^2+b^2)(d^2+e^2)=  (r^4+s^4)(u^4+w^4).
$$
But $ r^4,s^4 u^4,w^4$ cannot be legs of an integer Pythagorean triangle.
The outcome of the above argument is that there are no perfect cuboid among Euler bricks obtained using two primitive Pythagorean triangles.

\medskip

After completing this note, I found in Sierpinski's book \cite{Se} (pages 90--91 of the Russian translation)
the description the method of multiplication of Pythagorean triangles, which was put forward Ch.\,L.~Schedd \cite{Sch} in 1949.
Of course, I cannot prove that the idea was not borrowed by me, but I hope that the reader,
after having read the book, will know what it is all about. Nevertheless, Ch.\,L.~Schedd may be regarded as the founder of this method.


\begin{thebibliography}{9}
\bibitem{Dav} H.~Davenport, \textit{The Higher Arithmetic} (Cambridge, Cambridge University Press, 1999; Moscow, Nauka, 1965)
\bibitem{Ostrik}  V.\,V.~Ostrik, M.\,A.~Tsfasman, \textit{Algebraic Geometry and Number Theory: Rational and Elliptic Curves} (Moscow, MTsNMO, 2001).
\bibitem{Sch} Ch.\,L.~Schedd, ``A Hypothense Common to 64 Primitive Right Triangles,'' Scripta Math. {\bf 15}, 1949
\bibitem{Se} W. Sierpi\'nski,  \textit{Pythagorean Triangles} (Warszaw, State Scientific Publishers, 1954; New York, Academic Press, 1962; Moscow, Uchpedgiz, 1959).
\bibitem{Lee} J. Leech, ``The rational cuboid revisited,'' Amer. Math. Monthly {\bf 84}\,(7), 518--533 (1977).

\end{thebibliography}
\end{document}